\documentclass{amsart}
\usepackage{latexsym,amsmath,amssymb}

\theoremstyle{definition}
\theoremstyle{remark}

\numberwithin{equation}{section}

\begin{document}

\title[WEIGHTED COMPOSITION OPERATORS]
{Hyperexpansive WEIGHTED COMPOSITION OPERATORS}

\author{\sc y. estaremi}
\address{\sc  y. estaremi}
\email{estaremi@gmail.com}
\address{Department of Mathematics, Payame Noor University, p. o. box: 19395-3697, Tehran, Iran.}

\thanks{}

\thanks{}



\date{}

\dedicatory{}

\commby{}

\begin{abstract}
In this note unbounded hyperexpansive weighted
composition operators are investigated. As a consequence unbounded
hyperexpansive multiplication and composition operators are
characterized.

\noindent {}
\end{abstract}

\maketitle

\section{ \sc\bf Introduction and Preliminaries}

%
Weighted composition operators are a general class of operators and they appear
naturally in the study of surjective isometries on most of the function spaces, semigroup
theory, dynamical systems, Brennan’s conjecture, etc. This type of operators are a generalization of multiplication operators and
composition operators. The main subject in the study
of composition operators is to describe operator theoretic properties of $C_{\phi}$ in terms of function
theoretic properties of $\phi$. The book \cite{co} is a good reference for the theory of composition
operators. Weighted composition operators had
been studied extensively in past decades. The basic properties of weighted composition
operators on measurable function spaces are studied by Lambert
\cite{lambe, lamber}, Singh and Manhas \cite{sin}, Takagi
\cite{ta}, Hudzik and Krbec \cite{h},  Cui, Hudzik, Kumar and  Maligranda \cite{c}, Arora \cite{ar}, Piotr Budzynski, Zenon Jan Jablonski, Il Bong Jung and Jan Stochel \cite{bjjs} and some
other mathematicians.\\

In this paper we consider unbounded weighted composition operators on the Hilbert space $L^2(\Sigma)$ and study hyperexpansive weighted
composition operators. As a consequence hyperexpansive multiplication and composition operators are
characterized.

Let $\mathcal{H}$ be stand for a Hilbert space and $B(\mathcal{H})$ for the Banach algebra of all bounded operators
on $\mathcal{H}$. By an operator on $\mathcal{H}$ we understand a
linear mapping $T:\mathcal{D}(T)\subseteq \mathcal{H}\rightarrow
\mathcal{H}$ defined on a linear subspace $\mathcal{D}(T)$ of
$\mathcal{H}$ which is called the domain of $T$. Set
$\mathcal{D}^{\infty}(T)=\cap^{\infty}_{n=1}\mathcal{D}(T^n)$.
Given an operator $T$ on $\mathcal{H}$, we define the graph norm
$\|.\|_{T}$ on $\mathcal{D}(T)$ by
$$\|f\|^2_{T}=\|f\|^2+\|Tf\|^2, \ \ \ \ f\in \mathcal{D}(T).$$
The next proposition can be easily deduced from the closed graph
theorem.\\

\vspace*{0.3cm} {\bf Proposition 1.1.} If $T$ is a closed operator
on $\mathcal{H}$ such that $T(\mathcal{D}(T)\subseteq
(\mathcal{D}(T)$, then $T$ is a bounded operator on the Hilbert
space $(\mathcal{D}(T), \|.\|_{T})$.\\

For an operator $T$ on $\mathcal{H}$ we set

$$\Theta_{T,n}(f)=\sum_{0\leq i\leq n}(-1)^i\left(%
\begin{array}{c}
  n \\
  i \\
\end{array}%
\right)\|T^i(f)\|^2, \ \ \ \ \ f\in \mathcal{D}(T^n), \ \
n\geq1.$$
We recall that an operator $T$ on $\mathcal{H}$ is:\\

(i) $k$-isometry $(k\geq1)$ if $\Theta_{T,k}(f)=0$ for $f\in
\mathcal{D}(T^k)$,\\

(ii) $k$-expansive $(k\geq1)$ if $\Theta_{T,k}(f)\leq0$ for $f\in
\mathcal{D}(T^k)$,\\

(iii) $k$-hyperexpansive $(k\geq1)$ if $\Theta_{T,n}(f)\leq0$ for
$f\in
\mathcal{D}(T^n)$ and $n=1,2,...,k$.\\

(iv) completely hyperexpansive  if $\Theta_{T,n}(f)\leq0$ for
$f\in
\mathcal{D}(T^n)$ and $n\geq1$.\\

\section{ \sc\bf Hyperexpansive weighted composition operators  }

Let $(X,\Sigma,\mu)$ be a $\sigma$-finite measure space. We denote the
collection of (equivalence classes modulo sets of zero measure of)
$\Sigma$-measurable complex-valued functions on $X$ by
$L^0(\Sigma)$ and the support of a function $f\in L^0(\Sigma)$ is
defined as $S(f)=\{x\in X; f(x)\neq 0\}$. We also
adopt the convention that all comparisons between two functions or
two sets are to be interpreted as holding up to a $\mu$-null set.
Denote by $L^2(\mu)$ the Hilbert space of all square summable (with respect to $\mu$) $\Sigma$-measurable complex functions on $X$.\\
For each $\sigma$-finite subalgebra $\mathcal{A}$ of $\Sigma$, the
conditional expectation, $E^{\mathcal{A}}(f)$, of $f$ with respect
to $\mathcal{A}$ is defined whenever $f\geq0$ almost everywhere or
$f\in L^2$.  For a sub-$\sigma$-finite algebra
$\mathcal{A}\subseteq\Sigma$, the conditional expectation
operator associated with $\mathcal{A}$ is the mapping
$f\rightarrow E^{\mathcal{A}}f$, defined for all non-negative $f$
as well as for all $f\in L^2(\Sigma)$,
where $E^{\mathcal{A}}f$, by the Radon-Nikodym theorem, is the
unique $\mathcal{A}$-measurable function satisfying
$$\int_{A}fd\mu=\int_{A}E^{\mathcal{A}}fd\mu, \ \ \ \forall A\in \mathcal{A} .$$
 As an operator on
$L^{2}({\Sigma})$, $E^{\mathcal{A}}$ is an idempotent and
$E^{\mathcal{A}}(L^2(\Sigma))=L^2(\mathcal{A})$. If there is no
possibility of confusion we write $E(f)$ in place of
$E^{\mathcal{A}}(f)$ \cite{rao,z}.\\

For a complex $\Sigma$-measurable function $u$ on $X$. Define
the measure $\mu_u:\Sigma\rightarrow [0, \infty]$  by

$$\mu_u(E)=\int_{E}|u|^2d\mu, \ \ \ \ E\in \Sigma.$$

It is clear that the measure $\mu_u$ is also $\sigma$-finite. By
the Radon-Nikodym theorem, if $\mu_{u}\circ\phi^{-1}\ll\mu$, then
there exists a unique (up to a.e. $\mu$ equivalence)
$\Sigma$-measurable function $J:X\rightarrow [0,\infty]$ such that

$$\mu_{u}(\phi^{-1}(E))=\mu_{u}\circ\phi^{-1}(E)=\int_{E}Jd\mu,  \ \ \ \ E\in \Sigma.$$
If $\mu\circ \phi^{-1}\ll\mu$, then $\mu_u\circ \phi^{-1}\ll\mu$. So, by definition of $\mu_{u}\circ\phi^{-1}$ and applying conditional
expectation with respect to $\phi^{-1}(\Sigma)$, we get that
$J=hE(|u|^2)\circ \phi^{-1}$, where $h$ is the Radon-Nykodim derivative $\frac{d\mu\circ\phi^{-1}}{d\mu}$.\\

Let $(X,\Sigma,\mu)$ be a $\sigma$-finite measure space, $u$ be a $\Sigma$-measurable complex function and suppose that $\phi$ is a
mapping from $X$ into $X$ which is measurable (i.e.
$\phi^{-1}(\Sigma)\subseteq \Sigma$).  Define the operator $uC_{\phi}:\mathcal{D}(uC_{\phi})\subseteq L^2(\mu)\rightarrow L^2(\mu)$ by

$$\mathcal{D}(uC_{\phi})=\{f\in
L^2(\mu):u.f\circ\phi \in L^2(\mu)\},$$
$$uC_{\phi}(f)=u.f\circ\phi.$$
Of course such operators may not be well-defined. One can see by direct computation that if $\mu_{u}\circ\phi^{-1}\ll\mu$, then $uC_{\phi}$ is well-defined. And so, if $\phi$ is a non-singular transformation, then the operator $uC_{\phi}$ is well-defined. Well-defined operators of the form $uC_{\phi}(f)=u.f\circ\phi$ acting in
$L^2(\mu)=L^2(X,\Sigma,\mu)$ with $\mathcal{D}(uC_{\phi})=\{f\in
L^2(\mu):u.f\circ\phi \in L^2(\mu)\}$ are called weighted
composition operators. If $\mu\circ \phi^{-1}\ll\mu$, then for every $f\in \mathcal{D}(uC_{\phi})$ we have
\begin{align*}
\|uC_{\phi}(f)\|^2&=\int_{X}|u|^2|f\circ
\phi|^2d\mu\\
&=\int_{X}E(|u|^2)|f|^2\circ
\phi d\mu\\
&=\int_{X}hE(|u|^2)\circ \phi^{-1}|f|^2d\mu.
\end{align*}

By induction we get that for every $n\geq 1$

\begin{align*}
\|(uC_{\phi})^n(f)\|^2&=\int_{X}|u_{\phi,n}|^2|f\circ
\phi^n|^2d\mu\\
&=\int_{X}J_{n}|f|^2d\mu,
\end{align*}

for all $f\in\mathcal{D}((uC_{\phi})^n)$.
Where $u_{\phi,n}=u.u\circ\phi.u\circ\phi^2...u\circ\phi^{n-1}$, $J_n=hE(J_{n-1}|u|^2)\circ \phi^{-1}$, $h$ is the Radon-Nykodim derivative $\frac{d\mu\circ \phi^{-1}}{\mu}$, $E$ is conditional expectation with respect to $\phi^{-1}(\Sigma)$ and $J_0=1$.\\

\vspace*{0.3cm} {\bf Lemma 2.1.} Let $w=1+J$ and $d\nu=wd\mu$. Then we have \\

(a) $S(w)=X$ and $L^2(\nu)=\mathcal{D}(uC_{\phi})$,\\

(b) And also, the followings are equivalent;\\

(i) $uC_{\phi}$ is densely defined.\\

(ii) $J<\infty$ a.e. $\mu$.\\

\vspace*{0.3cm} {\bf Proof.} (a) Let $f$ be a measurable function on $X$. We have

\begin{align*}
\|f\|^2_{\mu}+\|uf\circ \phi\|^2&=\int_{X}|f|^2d\mu+\int_{X}|uf\circ \phi|^2d\mu\\
&= \int_{X}(1+J)|f|^2d\mu=\|f\|^2_{\nu}.
\end{align*}

This means that, $f\in \mathcal{D}(uC_{\phi})$ if and only if $f\in L^2(\nu)$. So $L^2(\nu)=\mathcal{D}(uC_{\phi})$.\\

(b) $(i)\rightarrow (ii)$ Set $F=\{J=\infty\}$. By  (a), $f\mid_{F}=0$ a.e. $\mu$ for every $f\in\mathcal{D}(uC_{\phi})$. This and (i) implies that  $f\mid_{F}=0$ a.e. $\mu$ for every $f\in L^2(\mu)$. So we have $\chi_{A\cap F}=0$ a.e. $\mu$ for all $A\in \Sigma$ with $\mu(A)<\infty$. By the $\sigma$-finiteness of $\mu$ we have $\chi_{F}=0$ a.e. i.e $\mu(E)=0$.\\

$(ii)\rightarrow (i)$ Here we prove that $L^2(\nu)$ is dense in $L^2(\mu)$. Suppose that $f\in L^2(\mu)$ such that $\langle f,g \rangle=\int_{X}f.\bar{g}d\mu=0$ for all $g\in L^2(\nu)$. For $A\in \Sigma$ we set $A_n=\{x\in A:w(x)\leq n\}$. It is clear that $A_n\subseteq A_{n+1}$ and $X=\cup^{\infty}_{n=1}A_n$. Since $(X,\Sigma,\mu)$ is $\sigma$-finite, hence $X=\cup^{\infty}_{n=1}X_n$ with $\mu(X_n)<\infty$. If we set $B_n=A_n\cap X_n$, then $B_n\nearrow A$ and so $f.\chi_{B_n}\nearrow f.\chi_{A}$ a.e. $\mu$. Since $\nu(B_n)\leq (n+1) \mu(B_n)<\infty$, we have $\chi_{B_n}\in L^2(\nu)$ and by our assumption $\int_{B_n}fd\mu=0$. Therefore by Fatou's lemma we get that $\int_{A}fd\mu=0$. Thus for all $A\in \Sigma$ we have $\int_{A}fd\mu=0$. This means that $f=0$ a.e. $\mu$ and so $L^2(\nu)$ is dense in $L^2(\mu)$.\\

%

If all functions $J_i=hE(J_{i-1}|u|^2)\circ \phi^{-1}, i=1,. .
. ,n$, are finite valued, where $h_i$ is the Radon-Nykodim
derivative $\frac{d\mu\circ\phi^{-i}}{d\mu}$, then we set
 $$\triangle_{J,n}(x)=\sum_{0\leq i\leq n}(-1)^i
\left(%
\begin{array}{c}
  n \\
  i \\
\end{array}%
\right)J_i(x).
$$\\

\vspace*{0.3cm} {\bf Proposition 2.2.} If $\mathcal{D}(uC_{\phi})$
is dense in $L^2(\Sigma)$, then the following conditions are
equivalent:\\

(i) $uC_{\phi}(\mathcal{D}(uC_{\phi}))\subseteq
\mathcal{D}(uC_{\phi})$.\\

(ii) There exists $c>0$ such that $J_2\leq c(1+J_1)$ a.e. $\mu$.\\

\vspace*{0.3cm} {\bf Proof.} $(i)\rightarrow (ii)$. Since
$uC_{\phi}$ is closed, densely defined and
$uC_{\phi}(\mathcal{D}(uC_{\phi}))\subseteq
\mathcal{D}(uC_{\phi})$, then by closed graph theorem $uC_{\phi}$
is a bounded operator on
$(\mathcal{D}(uC_{\phi}),\|.\|_{uC_{\phi}})$.  Hence there exists
$c>0$ such that $\|uC_{\phi}(f)\|^2_{uC_{\phi}}\leq
c\|f\|^2_{uC_{\phi}}$ for $f\in \mathcal{D}(uC_{\phi})$. By
replacing $f$ with $uC_{\phi}(f)$ we have
\begin{align*}
\|(uC_{\phi})^2(f)\|^2&\leq \|uC_{\phi}(f)\|^2+\|(uC_{\phi})^2(f)\|^2\\
&\leq c(\|f\|^2+\|uC_{\phi}(f)\|^2)
\end{align*}

i.e,
\begin{align*}
\int_{X}J_2|f|^2d\mu&\leq c(\int_{X}|f|^2d\mu+\int_{X}J_1|f|^2d\mu)\\
&=\int_{X}c(1+J_1)|f|^2d\mu.
\end{align*}
This implies that for all $f\in \mathcal{D}(uC_{\phi})$ and also
for all $f\in \overline{\mathcal{D}(uC_{\phi})}=L^2(\mu)$ we have

$$\int_{X}(c(1+J_1)-J_2)|f|^2d\mu\geq0$$
and so $J_2\leq c(1+J_1)$ a.e. $\mu$.\\

$(ii)\rightarrow (i)$. Let $f\in\mathcal{D}(uC_{\phi})$. Then by
assumption $J_2\leq c(1+J_1)$ a.e. $\mu$, we have
\begin{align*}
\int_{X}|(uC_{\phi})^2(f)|^2d\mu&=\int_{X}J_2|f|^2d\mu\\
&\leq
c(\int_{X}|f|^2d\mu
+\int_{X}J_1|f|^2d\mu)\\
&=c(\|f\|^2+\|uC_{\phi}(f)\|^2)<\infty.
\end{align*}
Therefore $uC_{\phi}(f)\in \mathcal{D}(uC_{\phi})$.\\

 \vspace*{0.3cm} {\bf Remark 2.3.} If $uC_{\phi}(\mathcal{D}(uC_{\phi}))\subseteq
\mathcal{D}(uC_{\phi})$ and $d\nu=(1+J_1)d\mu$, then
$(X,\Sigma,\nu)$ is a $\sigma$-finite measure space, $\nu\circ
\phi^{-1}$ is absolutely continuous with respect to $\nu$,
$L^2(\nu)=\mathcal{D}(uC_{\phi})$, $\|.\|_{L^2(\nu)}$ is the graph
norm of $uC_{\phi}$ (considered as an operator in $L^2(\mu)$), and
$uC_{\phi}$ is a bounded weighted composition operator acting on
$L^2(\nu)$. Furthermore, if $uC_{\phi}$ is $k$-isometric ( resp.
$k$-expansive, $k$-hyperexpansive), then so is $uC_{\phi}$ as an
operator on $L^2(\nu)$.\\

If all functions $u^{2i}$ and $h_i$ for $i=1,. . . ,n$ are finite
valued, then we set
 $$\triangle_{u,n}(x)=\sum_{0\leq i\leq n}(-1)^i
\left(%
\begin{array}{c}
  n \\
  i \\
\end{array}%
\right)u^{2i}(x),
$$

$$\triangle_{h,n}(x)=\sum_{0\leq i\leq n}(-1)^i
\left(%
\begin{array}{c}
  n \\
  i \\
\end{array}%
\right)h_i(x).
$$\\

\vspace*{0.3cm} {\bf Corollary 2.4.} If $\mathcal{D}(C_{\phi})$ is
dense in $L^2(\Sigma)$, then the following conditions are
equivalent:\\

(i) $C_{\phi}(\mathcal{D}(C_{\phi}))\subseteq
\mathcal{D}(C_{\phi})$.\\

(ii) There exists $c>0$ such that $h_2\leq c(1+h_1)$ a.e. $\mu$.\\

\vspace*{0.3cm} {\bf Corollary 2.5.} If $\mathcal{D}(M_u)$ is dense in
$L^2(\Sigma)$, then the following conditions are
equivalent:\\

(i) $M_u(\mathcal{D}(M_u))\subseteq
\mathcal{D}(M_u)$.\\

(ii) There exists $c>0$ such that $u^4\leq c(1+u^2)$ a.e. $\mu$.\\

\vspace*{0.3cm} {\bf Proposition 2.6.} If
$\mathcal{D}((uC_{\phi})^n)$
is dense in $L^2(\mu)$ for a fixed $n\geq1$, then:\\

(i) $uC_{\phi}$ is $k$-expansive if and only if
$\triangle_{J,n}(x)\leq 0$ a.e. $\mu$.\\

(ii) $uC_{\phi}$ is $k$-isometry if and only $\triangle_{J,n}(x)=0$ a.e. $\mu$.\\

\vspace*{0.3cm} {\bf Proof.}(i). Since
$\|(uC_{\phi})^i(f)\|^2=\int_{X}J_i|f|^2d\mu$ for all $f\in
\mathcal{D}((uC_{\phi})^i)$, we have
\begin{align*}
\sum_{0\leq i\leq n}(-1)^i
\left(%
\begin{array}{c}
  n \\
  i \\
\end{array}%
\right)\|(uC_{\phi})^i(f)\|^2&=\sum_{0\leq i\leq n}(-1)^i
\left(%
\begin{array}{c}
  n \\
  i \\
\end{array}%
\right)\int_{X}J_i|f|^2d\mu\\
&=\int_{X}\left(\sum_{0\leq i\leq
n}(-1)^i
\left(%
\begin{array}{c}
  n \\
  i \\
\end{array}%
\right)J_i\right)|f|^2d\mu\\
&=\int_{X}\triangle_{J,n}(x)|f|^2d\mu,
\end{align*}

for all $f\in \mathcal{D}((uC_{\phi})^n)$. Since $(uC_{\phi})^n$ is densely defined, then we get that
$uC_{\phi}$ is $k$-expansive if and only if
$\triangle_{J,n}(x)\leq 0$ a.e. $\mu$.\\

(ii) Likewise we have $uC_{\phi}$ is $k$-isometry if and only $\triangle_{J,n}(x)=0$ a.e. $\mu$.\\

\vspace*{0.3cm} {\bf Corollary 2.7 } If $\mathcal{D}((C_{\phi})^n)$
is dense in $L^2(\mu)$ for a fixed $n\geq1$, then:\\

(i) $C_{\phi}$ is $k$-expansive if and only if
$\triangle_{h,n}(x)\leq 0$ a.e. $\mu$.\\

(ii) $C_{\phi}$ is $k$-isometry if and only $\triangle_{h,n}(x)=0$ a.e. $\mu$.\\

\vspace*{0.3cm} {\bf Corollary 2.8.} If $\mathcal{D}((M_u)^n)$
is dense in $L^2(\mu)$ for a fixed $n\geq1$, then:\\

(i) $M_u$ is $k$-expansive if and only if
$\triangle_{u,n}(x)\leq 0$ a.e. $\mu$.\\

(ii) $M_u$ is $k$-isometry if and only $\triangle_{u,n}(x)=0$ a.e. $\mu$.\\

\vspace*{0.3cm} {\bf Proposition 2.9.} If
$\mathcal{D}((uC_{\phi})^2)$ is dense in $L^2(\mu)$ and
$uC_{\phi}$ is 2-expansive, then:\\

(i) $uC_{\phi}$ leaves its domain invariant:\\

(ii) $J_k\geq J_{k-1}$ a.e. $\mu$ for all $k\geq1$.\\

\vspace*{0.3cm} {\bf Proof.} (i). By the Proposition 2.3 we get
that $J_2\leq 2J_1-1$. Hence for every $f\in
\mathcal{D}(uC_{\phi})$ we have
\begin{align*}
\|(uC_{\phi})^2(f)\|^2&=\int_{X}J_2|f|^2d\mu\\
&\leq 2\int_{X}J_1|f|^2d\mu-\int_{X}|f|^2d\mu<\infty,
\end{align*}
so $uC_{\phi}(f)\in \mathcal{D}(uC_{\phi})$.\\

(ii) Since $uC_{\phi}$ leaves its domain invariant, then
$\mathcal{D}(uC_{\phi})\subseteq \mathcal{D}^{\infty}(uC_{\phi})$.
So by lemma 3.2 (iii) of \cite{jjst} we get that
$\|(uC_{\phi})^k(f)\|^2\geq \|(uC_{\phi})^{k-1}(f)\|^2$ for all
$f\in \mathcal{D}(uC_{\phi})$ and $k\geq1$ we have
$$\int_{X}(J_k-J_{k-1})|f|^2d\mu\geq 0, \ \ \ \ \ \ f\in \mathcal{D}(uC_{\phi}),$$
so this leads to $J_k\geq J_{k-1}$ a.e. $\mu$.\\

\vspace*{0.3cm} {\bf Corollary 2.10.} If $\mathcal{D}((C_{\phi})^2)$ is
dense in $L^2(\mu)$ and
$C_{\phi}$ is 2-expansive, then:\\

(i) $C_{\phi}$ leaves its domain invariant:\\

(ii) $h_k\geq h_{k-1}$ a.e. $\mu$ for all $k\geq1$.\\

\vspace*{0.3cm} {\bf Corollary 2.11.} If $\mathcal{D}((M_u)^2)$ is
dense in $L^2(\mu)$ and
$M_u$ is 2-expansive, then:\\

(i) $M_u$ leaves its domain invariant:\\

(ii) $u^{2k}\geq u^{2(k-1)}$ a.e. $\mu$ for all $k\geq1$.\\

Recall that a real-valued map $\varphi$ on $\mathbb{N}$ is said to
be completely alternating if $\sum_{0\leq i\leq n}(-1)^i
\left(%
\begin{array}{c}
  n \\
  i \\
\end{array}%
\right)\varphi(m+i)\leq0$ for all $m\geq 0$ and $n\geq 1$. The
next theorem is a direct consequence of proposition 2.3 and 2.4.\\

\vspace*{0.3cm} {\bf Theorem 2.12.} If $\mathcal{D}((uC_{\phi})^2)$
is dense in $L^2(\mu)$ and $k\geq1$
is fixed, then:\\

(i) $uC_{\phi}$ is $k$-hyperexpansive if and only if
$\triangle_{J,i}(x)\leq 0$ a.e. $\mu$ for $i=1,...,k$.\\

(ii) $uC_{\phi}$ is completely hyperexpansive if and only if
$\{J_i\}^{\infty}_{i=0}$ is a completely alternating sequence for
almost every $x\in X$.\\

\vspace*{0.3cm} {\bf Corollary 2.13.} If $\mathcal{D}((C_{\phi})^2)$ is
dense in $L^2(\mu)$ and $k\geq1$
is fixed, then:\\

(i) $C_{\phi}$ is $k$-hyperexpansive if and only if
$\triangle_{h,i}(x)\leq 0$ a.e. $\mu$ for $i=1,...,k$.\\

(ii) $C_{\phi}$ is completely hyperexpansive if and only if
$\{h_i\}^{\infty}_{i=0}$ is a completely alternating sequence for
almost every $x\in X$.\\

\vspace*{0.3cm} {\bf Corollary 2.14.} If $\mathcal{D}((M_u)^2)$ is
dense in $L^2(\mu)$ and $k\geq1$
is fixed, then:\\

(i) $M_u$ is $k$-hyperexpansive if and only if
$\triangle_{u,i}(x)\leq 0$ a.e. $\mu$ for $i=1,...,k$.\\

(ii) $M_u$ is completely hyperexpansive if and only if
$\{u^{2i}\}^{\infty}_{i=0}$ is a completely alternating sequence
for
almost every $x\in X$.\\

Notice that in the same way we can characterize alternatingly
hyperexpansive weighted composition operators.\\

We say that the $\sigma$-algebra $\phi^{-1}(\Sigma)$ is
essentially all of $\Sigma$ with respect to $\mu_u$ if and only of
given $A\in \Sigma$ there is $B\in \Sigma$ with the symmetric
difference $\phi^{-1}(B)\bigtriangleup A=(\phi^{-1}(B)\setminus
A)\cup (A\setminus\phi^{-1}(B))$ having
$\mu_u(\phi^{-1}(B)\bigtriangleup A)=0$.\\

 The following
proposition characterizes 2-expansive weighted composition
operators on the measure space $(X,\Sigma,\mu)$ such that
$\mu_u(X)<\infty$.\\

\vspace*{0.3cm} {\bf Theorem 2.15.} Let $uC_{\phi}$ be
2-expansive operator.\\

(i) Let $(X,\Sigma,\mu)$ is an infinite measure space such
that $\mu_u(X)<\infty$ and $\mathcal{D}((uC_{\phi})^2)$ is dense
in $L^2(\mu)$.\\

 (ii) Let $(X,\Sigma,\mu)$ is a measure space
such that $\mu_u(X)<\infty$, $u\leq1$ a.e. $\mu$ and $\mathcal{D}((uC_{\phi})^2)$ is dense in $L^2(\mu)$.\\

If the conditions (i) or (ii) holds, then $uC_{\phi}$ is an isometry.\\

(iii) If $uC_{\phi}$ is densely defined, $u\neq0$ a.e. $\mu$ and
the sigma algebra $\phi^{-1}(\Sigma)$ is essentially all of
$\Sigma$, with respect
to $\mu$, then $uC_{\phi}$ is a unitary operator.\\

\vspace*{0.3cm} {\bf Proof.} (i) It follows from Proposition 2.4
that $uC_{\phi}$ leaves its domain invariant and $J_1\geq1$ a.e.
$\mu$. Suppose that (a) holds and contrary to our claim, there
exists $B\subseteq X$ such that $\mu(B)>0$ and $J_1\geq
\epsilon+1$ on $B$ for some $\epsilon>0$. Then we have

$$\infty>\mu_u(X)=\mu_u(\phi^{-1}(X\setminus B))+\mu_u(\phi^{-1}(B))$$$$\geq
(1+\epsilon)\mu(B)+\mu(X\setminus B)>\mu(X),$$ which is a
contradiction. Thus $J_1=1$ a.e. $\mu$.\\

If (ii) holds, then by the same method we conclude that
$\mu(X)>\mu(X)$ which is a contradiction.  These imply that
$uC_{\phi}$ is an isometry.\\

(iii) Let $u\neq0$ a.e. $\mu$ and
the sigma algebra $\phi^{-1}(\Sigma)$ be essentially all of
$\Sigma$, with respect to $\mu$. This implies that $uC_{\phi}$ is dense range. Then by [\cite {jjst}, proposition 3.5] we get that $uC_{\phi}$ is unitary.\\

\vspace*{0.3cm} {\bf Corollary 2.16.} Let $C_{\phi}$ be
2-expansive operator.\\

(i) If $(X,\Sigma,\mu)$ is a finite measure space and
$\mathcal{D}((C_{\phi})^2)$ is dense
in $L^2(\mu)$, then $C_{\phi}$ is an isometry.\\

 (ii) If $C_{\phi}$ is densely defined and
the sigma algebra $\phi^{-1}(\Sigma)$ is essentially all of
$\Sigma$, then $C_{\phi}$ is a unitary operator.\\

\vspace*{0.3cm} {\bf Corollary 2.17.} Let $M_u$ be
2-expansive operator.\\

(i) Let $(X,\Sigma,\mu)$ is an infinite measure space such
that $\mu_u(X)<\infty$ and $\mathcal{D}((M_u)^2)$ is dense
in $L^2(\mu)$.\\

 (ii) Let $(X,\Sigma,\mu)$ is a measure space
such that $\mu_u(X)<\infty$, $u\leq1$ a.e. $\mu$ and $\mathcal{D}((M_u)^2)$ is dense in $L^2(\mu)$.\\

If the conditions (i) or (ii) holds, then $M_u$ is an isometry.\\

(iii) If $M_u$ is densely defined, $u\neq0$ a.e. $\mu$, then $M_u$ is a unitary operator.\\

Now, let $m=\{m_n\}_{n=1}^{\infty}$
be a sequence of positive real numbers. Consider the space
$l^2(m)=L^2(\mathbb{N}, 2^{\mathbb{N}},\mu)$, where
$2^{\mathbb{N}}$ is the power set of natural numbers and $\mu$ is
a measure on $2^{\mathbb{N}}$ defined by $\mu(\{n\})=m_n$. Let
$u=\{u_n\}_{n=1}^\infty$ be a sequence of complex numbers. Let $\varphi:\mathbb{N}\rightarrow \mathbb{N}$ be a non-singular measurable transformation; i.e. $\mu\circ
\varphi^{-1}\ll\mu$. Direct computation shows that
$$ h(k)=\frac{1}{m_k}\sum_{j\in
\varphi^{-1}({k})}m_j\ \ ,\ \ E_\varphi(f)(k)=\frac{\sum_{j\in
\varphi^{-1}(\varphi (k))}f_jm_j}{\sum_{j\in
\varphi^{-1}(\varphi(k))}m_j}\ , $$ for all non-negative sequence
$f=\{f_n\}_{n=1}^\infty$ and $k\in\mathbb{N}$.
So
$$J(k)=\frac{1}{m_k}\sum_{j\in
\varphi^{-1}(k)}|u_j|^2m_j.$$\\
This observations lets us to consider the weighted composition operators on discrete measure space $(\mathbb{N},\mu, \Sigma)$. If $uC_{\phi}$ is a weighted composition operator on $l^2(m)$, then\\
$$
\mathcal{D}(uC_{\phi})=\{f=\{f_n\}\in l^2(m):\sum^{\infty}_{n=0}(\sum_{j\in
\varphi^{-1}(n)}|u_j|^2m_j)|f_n|^2<\infty\}$$
\begin{align*}
\|uf\circ\phi\|^2&=\int_{X}|uf\circ\phi|^2d\mu\\
&=\sum^{\infty}_{n=0}J(n)m_n|f_n|^2\\
&=\sum^{\infty}_{n=0}(\sum_{j\in
\varphi^{-1}(n)}|u_j|^2m_j)|f_n|^2.
\end{align*}


\end{document}